\documentclass[11pt,twoside]{amsart}
\usepackage{amsmath}
\usepackage{amssymb}
\usepackage{amsfonts}
\usepackage{amscd}

\newcommand\op{\operatorname}

\newcommand\Aut{\op{Aut}}

\newcommand\Del{\op{Del}}
\newcommand\depth{\op{depth}}

\newcommand\for{\op{for}}

\newcommand\GL{\op{GL}}
\newcommand\Hilb{\op{Hilb}}
\newcommand\Hom{\op{Hom}}
\newcommand\id{\op{id}}

\newcommand\LEV{\op{LEV}}

\newcommand\red{\op{red}}

\newcommand\RIG{\op{RIG}}

\newcommand\Spec{\op{Spec}\,}

\newcommand\tsr{\op{tsr}}

\newcommand{\wcL}{{\widetilde {\cal L}}}

\newcommand{\wG}{{\widetilde G}}

\newcommand{\wP}{{\widetilde P}}
\newcommand{\wQ}{{\widetilde Q}}
\newcommand{\wR}{{\widetilde R}}

\newcommand{\bG}{{\bold G}}

\newcommand{\bP}{{\bold P}}

\newcommand{\bR}{{\bold R}}
\newcommand{\bZ}{{\bold Z}}

\newcommand{\XR}{X_{\bold R}}

\newcommand{\cG}{{\cal G}}

\newcommand{\cL}{{\cal L}}

\newcommand{\cN}{{\cal N}}
\newcommand{\cO}{{\cal O}}

\newtheorem{thm}[subsection]{Theorem}

\newtheorem{lemma}[subsection]{Lemma}
\newtheorem{cor}[subsection]{Corollary}
\newtheorem{maintheorem}{Theorem}

\theoremstyle{definition}
\newtheorem{defn}[subsection]{Definition}

\newtheorem{rem}[subsection]{Remark}
\ifx\undefined\pf
\newcommand{\pf}{\proof}
\fi
\ifx\undefined\cal
\newcommand{\cal}{\mathcal}
\fi
\ifx\undefined\Bbb
\newcommand{\Bbb}{\mathbb}
\fi
\ifx\undefined\bold
\newcommand{\bold}{\mathbf}
\fi
\ifx\undefined\frak
\newcommand{\frak}{\mathfrak}
\fi

\begin{document}
\title[A compactification of the moduli scheme]
{A compactification\\
of the moduli scheme of abelian varieties}
\author{Iku Nakamura}
\address{Department of Mathematics, 
Hokkaido University,Sapporo, 060-0810, Japan} 
%\email{nakamura\@math.sci.hokudai.ac.jp}
\email{nakamura@math.sci.hokudai.ac.jp}
\thanks{Research was supported in part by  the Grant-in-aid 
(No. 12304001, No. 11874001) for Scientific Research, 
JSPS}
%\date{\today}
%
\date{July 24, 2001}
\subjclass{Algebraic Geometry,
14J10, 14K10, 14K25}
\keywords{Moduli, Abelian variety,  
Torically stable quasi-abelian scheme, Coarse moduli,
Heisenberg group, Level structure}
\begin{abstract} We construct a canonical
compactification $SQ^{toric}_{g,K}$ of the moduli  
of abelian varieties over $\bZ[\zeta_N, 1/N]$ where 
$\zeta_N$ is a primitive $N$-th root of unity. 
It is very similar to,  but slightly
different from the compactification $SQ_{g,K}$ in 
\cite{Nakamura99}. 
Any degenerate abelian scheme 
on the boundary of $SQ^{toric}_{g,K}$ is one of the  
(torically) stable quasi-abelian schemes 
introduced in \cite{AN99}, which is reduced and singular. 
In contrast with it, some of degenerate abelian schemes 
on the boundary of $SQ_{g,K}$ are nonreduced schemes. 
\end{abstract}
\maketitle

\section{Introduction.}
\label{sec:Introduction}
In the article \cite{Nakamura99} a 
canonical compactification $SQ_{g,K}$ 
of the moduli scheme $A_{g,K}$ of abelian varieties with level structure 
was constructed by applying the geometric invariant theory \cite{MFK94}. 
It is a compactification of $A_{g,K}$
by all Kempf-stable degenerate abelian schemes. 
However some of the Kempf-stable degenerate abelian 
schemes are {\it nonreduced} 
by contrast with Deligne-Mumford stable curves.  
\par
The purpose of this article is 
to construct another canonical compactification $SQ^{toric}_{g,K}$ 
of $A_{g,K}$ by certain singular {\it reduced} 
degenerate abelian schemes only
instead of Kempf-stable ones. 
The new compactification $SQ^{toric}_{g,K}$ 
is very similar to $SQ_{g,K}$. 
The compactifications are as functors the same  
if $g\leq 4$, and different if $g\geq 8$ (or maybe if $g\geq 5$).
An advantage of $SQ^{toric}_{g,K}$ is that 
the {\it reduced} degenerate abelian schemes 
on the boundary $SQ^{toric}_{g,K}\setminus A_{g,K}$
are much simpler than Kempf-stable ones.
\par 
Let $R$ be a complete discrete valuation ring and $k(\eta)$
the fraction field of $R$. 
Given an abelian variety $(G_{\eta},\cL_{\eta})$ 
over $k(\eta)$,
we have Faltings-Chai's degeneration data of it
by a finite base change if necessary. 
Then there are two natural choices 
of $R$-flat projective degenerating  
families $(P,\cL)$ and $(Q,\cL)$ of abelian varieties  
with generic fibre isomorphic to $(G_{\eta},\cL_{\eta})$ where
$(Q,\cL)$ is the most naive choice, while $(P,\cL)$ 
is the normalization of 
$(Q,\cL)$ after a certain finite base change such that the closed fibre
$P_0$ is reduced. 
We call the closed fibre $(P_0,\cL_0)$ of $(P,\cL)$ 
a torically stable quasi-abelian scheme (abbr. TSQAS) \cite{AN99}, while 
we call $(Q_0,\cL_0)$ a projectively stable quasi-abelian scheme 
 (abbr. PSQAS) \cite{Nakamura99}.
\par
Let $(K,e_K)$ be a finite symplectic abelian group. 
Since $K\simeq\oplus_{i=1}^g(\bZ/e_i\bZ)^{\oplus 2}$ 
for some positive integers $e_i$ such that $e_i|e_{i+1}$, we define  
$e_{\op{min}}(K)=e_1$ and $e_{\op{max}}(K)=e_g$. 
Let $N=e_{\op{max}}(K)$.
The Heisenberg group $G(K)$ is by definition a central
extension of $K$ by the group $\mu_{N}$ of all $N$-th roots of
unity.  
The classical level $K$-structures on abelian varieties are generalized as 
level $G(K)$-structures on TSQASes. 
\begin{maintheorem}\label{thm:moduli/Z intro}   
If $e_{\min}(K)\geq 3$,  the functor
of $g$-dimensional torically 
stable quasi-abelian schemes with level $G(K)$-structure 
over reduced base algebraic spaces is coarsely represented 
by a complete reduced separated algebraic space $SQ^{toric}_{g,K}$ over 
$\bZ[\zeta_N, 1/N]$.
\end{maintheorem}

We prove 
the theorem with the help of
\cite{Nakamura99} and Keel-Mori \cite{KM97}. 
Alexeev \cite{Alexeev99} treats similar problems 
in a different formulation. 
\par
Here is an outline of our article. 
In Section~\ref{sec:degenerating families} 
we recall from \cite{Nakamura99} 
a couple of basic facts about 
degenerating families of abelian varieties. 
In Section~\ref{sec:level G(K) structure}  
we define a rigid $G(K)$-structure and a level $G(K)$-structure on a
TSQAS $(P_0,\cL_0)$ or their family. 
In Section~\ref{sec:The stable reduction theorem} 
we recall from \cite{Nakamura99} the stable
reduction theorem  
for TSQASes with rigid $G(K)$-structure.  
In Sections~\ref{sec:The scheme parametrizing TSQASes},
\ref{sec:The fibres over U3} and \ref{sec:The geometric quotient} we
prove the above theorem. \par
\medskip
\noindent 
{\it Acknowledgement.}
The author would like to thank Takeshi Sugawara for stimulating
discussions.
%%%%%%%%%%
%%%%%%%%%%
%%%%%%%%%%
%%%%%%%%%%
\section{Degenerating families of abelian varieties}
\label{sec:degenerating families}
%\end{document}

The purpose of this section is 
to recall basic facts about degenerating families 
of abelian varieties.  
To minimize the article we try to keep the same notation as 
\cite{Nakamura99}. 

\subsection{Grothendieck's stable reduction}
\label{subsec:Grothendieck stable red}
Let $R$ be a complete discrete valuation ring, 
$I$ the maximal ideal of $R$  and $S=\Spec R$. 
Let  $\eta$ be the generic point of $S$, 
$k(\eta)$ the fraction field of $R$
and  $k(0)=R/I$ the residue field. \par
Suppose we are given a $g$-dimensional polarized abelian variety 
$(G_{\eta}, \cL_{\eta})$ over $k(\eta)$ such that 
$\cL_{\eta}$ is symmetric, ample and rigidified (that is trivial) 
along the unit section. 
Then by Grothendieck's stable reduction theorem
\cite{SGA7}
$(G_{\eta}, \cL_{\eta})$  can be extended 
to a polarized semiabelian $S$-scheme $(G, \cL)$ 
with $\cL$ a rigidified relatively ample invertible sheaf on $G$ 
as the connected N\'eron model of $G_{\eta}$ 
by taking a finite extension $K'$ of $k(\eta)$ if necessary. 
The closed fibre $G_0$ is a semiabelian scheme over $k(0)$, namely
an extension of an abelian variety $A_0$
by a split torus $T_0$.\par
{\it From now on we restrict ourselves 
to the totally degenerate case, that is, the case where $A_0$ is trivial}
because by \cite{Nakamura99} 
there is no essentially new difficulty 
when we consider the case $A_0$ is nontrivial. 
Hence we assume that $G_0$ is a split $k(0)$-torus.
Let $\lambda(\cL_{\eta}):G_{\eta}\to G^t_{\eta}$ 
be the polarization (epi)morphism. 
By the universal property 
of the (connected) N\'eron model $G^t$ of 
$G^t_{\eta}$ we have an epimorphism $\lambda:G\to G^t$ 
extending $\lambda(\cL_{\eta})$. 
Hence 
the closed fibre of $G^t$ is also a split $k(0)$-torus.\par
Let $S_n=\Spec R/I^{n+1}$ and $G_n=G\times_SS_n$.
Associated to $G$ and $\cL$ are the formal scheme 
$G_{\for}=\lim G_n$ and an invertible sheaf 
$\cL_{\for}=\lim \cL\otimes R/I^{n+1}$. 
By our assumption that $G_0$ is 
a $k(0)$-split torus, $G_n$ turns out to be a 
multiplicative group scheme for every $n$ 
by \cite[p.~7]{FC90}. 
Thus the scheme 
$G_{\for}$ is a formal split $S$-torus. Similarly
$G^t_{\for}$ is a formal split $S$-torus.
Let $X:=\Hom_{\bZ}(G_{\for},(\bG_{m,S})_{\for})$, 
$Y:=\Hom_{\bZ}(G^t_{\for},(\bG_{m,S})_{\for})$ and 
$\wG:=\Hom_{\bZ}(X,\bG_{m,S})$, $\wG^t=\Hom(Y,\bG_{m,S})$. 
Then $\wG$ (resp. $\wG^t$) algebraizes
$G_{\for}$ (resp. $G^t_{\for}$). 
The morphism $\lambda:G\to G^t$ 
induces an injective homomorphism $\phi:Y\to X$ and 
an algebraic epimorphism $\widetilde\lambda:\wG\to \wG^t$.
For simplicity we identify the injection $\phi : Y\to X$ with 
the inclusion $Y\subset X$.  
\subsection{Fourier expansions}
In the totally degenerate case 
$G_{\for}$ (resp. $\wG$) is a formal split $S$-torus (resp. a
split $S$-torus). 
We choose and fix the coordinate $w^x$ $(x\in X)$ of $\wG$ satisfying
$w^xw^y=w^{x+y}$ $(\forall x,y\in X)$.
Since $\cL_{\for}$ 
is trivial on $G_{\for}$, we have
\begin{align*}
\Gamma(G_{\eta},\cL_{\eta})=\Gamma(G,\cL)\otimes k(\eta)& 
\hookrightarrow \Gamma(G_{\for},\cL_{\for})\otimes k(\eta)
\hookrightarrow \prod_{x\in X} 
k(\eta)\cdot w^x.\end{align*}
Therefore, any element
$\theta\in\Gamma(G_{\eta},\cL_{\eta})$ 
can be written as a formal Fourier series
$\theta= \sum_{x\in X} \sigma_x(\theta)w^x$
with $\sigma_x(\theta)\in k(\eta)$, 
which converges $I$-adically.

\begin{thm}{$\op{[Faltings--Chai90]}$}
\label{thm:Fourier_series_max_deg}Let $k(\eta)^*=k(\eta)\setminus \{0\}$.
There exists a function $a:Y\to k(\eta)^*$ and 
a bimultiplicative  function $b:Y\times X\to k(\eta)^*$ 
with the following properties: \begin{enumerate}
\item\label{a,b} $b(y,x)=a(x+y)a(x)^{-1}a(y)^{-1}$, 
$a(0)=1$\ \ $(\forall x\in X, \forall y\in Y)$,
\item\label{b} $b(y,z)=b(z,y)=a(y+z)a(y)^{-1}a(z)^{-1}$ 
$\ \ (\forall y,z\in Y)$,
\item\label{item:positive definite} $b(y,y)\in I \ \ (\forall y\ne 0)$, 
and for 
every $n\ge 0$, $a(y)\in I^n$ for almost all $y\in Y$, 
\item\label{k(eta) module} 
$\Gamma(G_\eta,\cL_\eta)$ is 
identified with the $k(\eta)$ vector subspace of formal Fourier series 
$\theta$ that satisfy $\sigma_{x+y}(\theta)
=a(y)b(y,x)\sigma_x(\theta)$ and $\sigma_x(\theta)\in k(\eta)$ 
$(\forall x\in X, \forall y\in Y)$. 
\end{enumerate}
\end{thm}

\begin{defn} By taking a finite base change of $S$ if necessary
the functions $b$ and $a$ can be extended respectively 
to $X\times X$ and $X$ 
so that the previous relations between 
$b$ and $a$ 
are still true on $X\times X$. Let $R^*=R\setminus \{0\}$ 
and $k(0)^*=k(0)\setminus \{0\}$.
Then we define integer-valued functions
$A : X \to \bZ$, $B : X\times X\to \bZ$ 
and $\bar b(y,x)\in R^*$, $\bar a(y)\in R^*$ by
\begin{align*}
B(y,x)&=val_s(b(y,x)),\quad
dA(\alpha)(x)=B(\alpha,x)+r(x)/2,\\
A(x)&=val_s(a(x))=B(x,x)/2+r(x)/2,\\
b(y,x)&=\bar b(y,x)s^{B(y,x)},\quad 
a(x)=\bar a(x)s^{(B(x,x)+r(x))/2}
\end{align*}
for some $r\in \Hom_{\bZ}(X,\bZ)$, where $B$ is positive definite by 
Theorem \ref{thm:Fourier_series_max_deg}
(\ref{item:positive definite}).   
We set $a_0=\bar a \mod I$ and $b_0=\bar b\mod I$.  
Hence $a_0(x), b_0(y,x)\in k(0)^*$.
\end{defn}
\begin{defn}
\label{defn:Delaunay cell} Let $X$ be a lattice of rank $g$, 
$X_{\bR}=X\otimes\bR$,
and let $B:X\times X\to \bZ$ be 
a positive definite symmetric integral bilinear form, 
which determines a distance 
 $\|\ \|_B$ on $X_{\bR}$ by 
$\|x\|_B:=\sqrt{B(x,x)}$ $(x\in X_{\bR})$.
For any $\alpha\in X_{\bR}$ we say 
that $a\in X$ is $\alpha$-nearest if
\begin{displaymath}
\|a-\alpha\|_B= \min\{ \|b-\alpha\|_B; b\in X \}. \end{displaymath}

We define {\em a Delaunay cell\/ $\sigma$} 
 to be the closed convex hull of all lattice elements 
which are $\alpha$-nearest for some $\alpha\in\XR$. 
All the Delaunay cells constitute 
a locally finite decomposition of $\XR$ which we call the Delaunay 
decomposition $\Del_B$.
Let $\Del:=\Del_B$, 
and $\Del(c)$ the set of all the Delaunay cells containing $c\in X$. 
For $\sigma\in\Del(c)$, we define 
$C(c,\sigma)$ to be the cone spanned by edges of $\sigma$ at $c$. 
See \cite[p.~662]{Nakamura99}. \end{defn}

\begin{defn}\label{defn:defn of tilde R} 
We define
\begin{align*}
\wR:&=
\ R[a(x)w^x\vartheta; x\in X]
\simeq \ R[\xi_x\vartheta; x\in X],\\
\xi_x:& =s^{B(x,x)/2+r(x)/2}w^x,\\
\zeta_{x,c}:&=
s^{B(\alpha(\sigma),x)+r(x)/2}w^x\quad
(x+c\in C(c,\sigma))
\end{align*}
where $\wR$ is the graded algebra with 
$\deg(a(x)w^x\vartheta)=1$ and $\deg a=0$ for $a\in R$, while 
$\sigma\in\Del(c)$ is 
a Delaunay $g$-cell with $x+c\in C(c,\sigma)$.  
Let $\wQ:=\op{Proj}(\wR)$. 
We define an action $S_y$ on $\wQ$ by 
\begin{equation*}\label{Sy}
S_y^*(a(x)w^x\vartheta)=a(x+y)w^{x+y}\vartheta\quad
\op{for}\ y\in Y.
\end{equation*}
which induces a natural action on $\wP$, the normalization 
of $\wQ$, denoted by the same $S_y$.  
By $\wcL$ we denote $O_{\op{Proj}}(1)$ on $\wQ$ as well as
its pullback to $\wP$. 
\end{defn}

\begin{thm}\label{thm:construction of P}
\begin{enumerate}
\item The quotients $(\wP_{\for},\wcL_{\for})/Y$ and
$(\wQ_{\for},\wcL_{\for})/Y$ are   
flat projective formal $S$-schemes.
\item 
There exist flat projective $S$-schemes 
$(P,\cL)$ and $(Q,\cL)$ such that 
their formal completions $(P_{\for},\cL_{\for})$ and
$(Q_{\for},\cL_{\for})$
along the closed fibres are respectively isomorphic to 
$(\wP_{\for},\wcL_{\for})/Y$ and $(\wQ_{\for},\wcL_{\for})/Y$.
\item $P$ is the normalization of $Q$.
\end{enumerate}
\end{thm} 
\begin{pf} 
This follows from  \cite[III,~5.4.5]{EGA} 
and \cite{Mumford72}. See also \cite{Nakamura99}.
\qed\end{pf}

\begin{defn}\label{defn:Del P0} Let $\Del(P_0)$ be the Delaunay 
decomposition corresponding to $P_0$.
By taking a finite base change of $S$ if necessary, we may assume that 
$dA(\alpha(\sigma))\in\Hom(X,\bZ)$ for any Delaunay cell $\sigma\in 
\Del(P_0)$. By \cite{AN99} this implies that $P_0$ is reduced. 
We call the closed fibre $(P_0,\cL_0)$ of
$(P,\cL)$ a  
torically stable quasi-abelian scheme (abbr. a TSQAS) over $k(0):=R/I$.\par
In what follows 
we always assume that $dA(\alpha(\sigma))\in\Hom(X,\bZ)$ 
for any $\sigma\in\Del(P_0)$.
Hence $P_0$ is reduced. 
\end{defn}

We quote a few theorems from \cite{AN99} and \cite{Nakamura99}.
\begin{thm}\label{thm:local structure of fibre P} 
Let $(\wP_0,\wcL_0)$ be the 
closed fibre of $(\wP,\wcL)$ and 
$\bar\zeta_{x,c}:=\zeta_{x,c}\otimes k$ 
the restriction to $P_0$. Then 
\begin{enumerate}
\item
$\wP_0$ is covered 
with infinitely many affine schemes of finite type
$$U_0(c):=\Spec k[\bar\zeta_{x,c}; x\in X]\quad (c\in X).$$
\item $R_0(c):=k[\bar\zeta_{x,c}; x\in X]$ 
is a $k$-algebra of finite type. 
Let $x_i\in X$.  
If $x_i$'s belong to one and the same Delaunay cell (resp. otherwise), 
then 
\begin{equation*}
\bar\zeta_{x_1,c}\cdots \bar\zeta_{x_m,c}
=\bar\zeta_{x_1+\cdots +x_m,c}\quad ({\it resp.}\  0).
\end{equation*}
\end{enumerate}
\end{thm}

\begin{thm}\label{thm:stratification} 
Let $(P_0,\cL_0)$ be the closed fibre of $(P,\cL)$. 
Let $\sigma$ and $\tau$ be Delaunay cells in $\Del(P_0)$.
\begin{enumerate} 
\item For each $\sigma\in\Del(P_0)$ there exists 
a $G$-invariant subscheme $O(\sigma)$ of  
$P_0$ which is a torus of dimension $\dim\sigma$ 
over $k(0)$,
\item $\sigma\subset\tau$ iff $O(\sigma)$ is contained 
in $\overline{O(\tau)}$, the closure of $O(\tau)$ in $P_0$,
\item $P_0=\bigcup_{\sigma\in\Del(P_0)\mod Y} O(\sigma)$.
\end{enumerate}
\end{thm}

\begin{thm}\label{thm:dimension} Let $n>0$. Then
\begin{enumerate}
\item 
$h^0(P_0,\cL^n_0)=[X:Y]n^g$, $h^i(P_0,\cL^n_0)=0$ $(i>0)$, and 
$$\Gamma(P_0,\cL_0) = 
\left\{\sum_{x\in X}\ c(x)\xi_x;\begin{matrix} 
c(x+y)=b_0(y,x)a_0(y)c(x)\\
c(x)\in k,\forall x\in X,\forall y\in Y
\end{matrix}
\right\}.
$$
\item $\cL^n_0$ is very ample for $n\geq 2g+1$.
\end{enumerate}\end{thm}

%%%%%%%%%%%%%%%%%%
%%%%%%%%%%%%%%%%%%
\subsection{The group schemes $G$ and $G^{\sharp}$}
\label{subsec:group scheme G and Gsharp}We review 
\cite[3.12]{Nakamura99} to recall the notation.
By choosing a suitable base change of $S$ we assume
$dA(\alpha(\sigma))\in\Hom(X,\bZ)$ 
for any $\sigma\in \Del_B$.  Then 
$P_0$ is reduced. 
Then $G$ is realized as an open subscheme of $P$. In fact,  
for any Delaunay $g$-cell 
$\sigma\in\Del(0)$, there is an open smooth subscheme $G(\sigma)\subset P$ 
such that 
\begin{enumerate}
\item $G(\sigma)\simeq G$, $G(\sigma)_{\eta}=P_{\eta}$,
$G(\sigma)_0=O(\sigma)$, 
\item $G(\sigma)_{\for}$ is a formal $S$-torus.
\end{enumerate}

We define $G^{\sharp}=G^{\sharp}(\sigma):=\cup_{x\in (X/Y)}
S_x(G(\sigma))\subset P$. 
Then $G^{\sharp}$ is a group scheme over $S$ such that 
$G^{\sharp}_{\eta}=P_{\eta}$. It 
is uniquely determined by $P$, independent of the choice of $\sigma$.
See \cite[4.13]{Nakamura99} for the detail.

%%%%%%%%%%%%%%%%
%%%%%%%%%%%%%%%%
\subsection{The Heisenberg group}
\label{subsec:Heisenberg group over k(eta)} 
Let $K(\cL_{\eta})$ be 
the kernel of $\lambda(\cL_{\eta}):G_{\eta}\to G^t_{\eta}$. It is 
a subgroup scheme of $G_{\eta}$ 
representing the functor defined by
\begin{equation*}\label{eq:K(cLeta)}
K(\cL_{\eta})(U)=\{x\in G_{\eta}(U);  
T_x^*(\cL_{\eta,U})\simeq \cL_{\eta,U}\} 
\end{equation*}
for an $k(\eta)$-scheme $U$ where $\cL_{\eta,U}$ 
is the pullback of $\cL_{\eta}$ to 
$(G_{\eta})\times_{k(\eta)} U$. 
Let $\cL^{\tsr}_{\eta}$ be the $\bG_m$-torsor on $G_{\eta}$
associated with the invertible sheaf $\cL_{\eta}$. 
Let $\cG(\cL_{\eta}):=(\cL^{\tsr}_{\eta})_{|K(\cL_{\eta})}$\ 
(which we call the Heisenberg group of $\cL_{\eta}$). 
Then $\cG(\cL_{\eta})$ is a group $k(\eta)$-scheme which 
represents the functor 
\begin{equation*}
\begin{aligned}\label{eq:calG(eLeta)}
\cG(\cL_{\eta})(U)&=\{(a,\phi);a\in K(\cL_{\eta})(U)\ \op{and}\\
&\hskip 1cm \phi\ 
\op{is\ an\ isomorphism\ of}\  
T^*_a(\cL_{\eta,U})\ \op{with}\ \cL_{\eta,U}\}
\end{aligned}
\end{equation*}
  
It is a central extension of 
$K(\cL_{\eta})$ by the $k(\eta)$-split torus 
$\bG_{m,\eta}$. We define the 
commutator form $e^{\cL_{\eta}}$ of $\cG(\cL_{\eta})$ by
\begin{equation*}\label{eq:eLeta}
e^{\cL_{\eta}}(\bar g,\bar h)=[g,h]:=ghg^{-1}h^{-1},\quad 
\op{for}\ \forall g,h\in\cG(\cL_{\eta}) 
\end{equation*}
where $\bar g$, $\bar h$ are the images of 
$g$ and $h$ in $K(\cL_{\eta})$.  It 
is bimultiplicative nondegenerate and alternating on $K(\cL_{\eta})$.
By \cite[V, 2.5.5]{MB85} (See also \cite[\S 23]{Mumford74}), 
$\Gamma(G_{\eta},\cL_{\eta})$ is 
an irreducible $\cG(\cL_{\eta})$-module of weight one, 
unique up to isomorphism 
by taking a finite extension of $k(\eta)$ if necessary.

\begin{lemma}\label{lemma:finite flat closure} The flat closure 
$K^{\sharp}_S(\cL_{\eta})$ of $K(\cL_{\eta})$ in $G^{\sharp}$ 
is finite over $S$. 
\end{lemma}

See \cite[Lemma~4.14]{Nakamura99}.

\begin{defn}\label{defn:GsharpLeta}
Let $\cL^{\tsr}$ be the $\bG_m$-torsor on $G^{\sharp}$
associated with the invertible sheaf $\cL_{|G^{\sharp}}$. 
Let $\cG^{\sharp}_S(\cL_{\eta}):
=\cL^{\tsr}_{|K^{\sharp}_S(\cL_{\eta})}$ and 
$e_S^{\sharp}$ an extension of $e^{\cL_{\eta}}$ to
$K^{\sharp}_S(\cL_{\eta})$. 
By \cite[IV, 7.1 (ii)]{MB85}
$\cG^{\sharp}_S(\cL_{\eta})$ is a group scheme over $S$ extending 
$\cG(\cL_{\eta})$, which is a central extension of
$K^{\sharp}_S(\cL_{\eta})$ 
by $\bG_{m,S}$ with $e_S^{\sharp}$ 
the commutator form.  The bimultiplicative form $e_S^{\sharp}$ 
on $K^{\sharp}_S(\cL_{\eta})$ is 
nondegenerate alternating by
\cite[IV, 2.4]{MB85} and 
by Lemma~\ref{lemma:finite flat closure}. 
\end{defn}
%%%%%%%%
%%%%%%%%
%%%%%%%%

\begin{defn}\label{defn:cal G(P,L)} 
We define 
\begin{align*}
K(P,\cL)&:=K^{\sharp}_S(\cL_{\eta}),\
\cG(P,\cL):=\cG^{\sharp}_S(\cL_{\eta}),\\
\cG(G_{\eta},\cL_{\eta})&:=\cG(\cL_{\eta})=
\cG(P,\cL)\otimes k(\eta),\\
K(P_0,\cL_0)&:=K(P,\cL)\otimes k(0),\
\cG(P_0,\cL_0):=\cG(P,\cL)\otimes k(0).
\end{align*}

The natural projection from 
$\cL^{\tsr}$ to $G^{\sharp}$ makes $\cG(P,\cL)$ 
a central extension  
of $K(P,\cL)$ by $\bG_{m,S}$ with commutator form $e^{\sharp}_S$
\begin{equation*}
1\to \bG_{m,S}\to \cG(P,\cL)\to K(P,\cL)\to 0.
\end{equation*}
We call $\cG(P,\cL)$ (resp. $\cG(P_0,\cL_0)$) 
the Heisenberg group scheme of $(P,\cL)$ (resp. $(P_0,\cL_0)$). 
Later in Definition~\ref{defn:finite Heisenberg group} 
we define a finite version $G(P,\cL)$ of $\cG(P,\cL)$. 
\end{defn}

\begin{lemma}\label{lemma:comparison of Gamma} 
Let $G^{\sharp}\subset P$ be the group $S$-scheme 
in \ref{subsec:group scheme G and Gsharp}. Then
\begin{enumerate}
\item
$\Gamma(Q,\cL)=\Gamma(P,\cL)=\Gamma(G^{\sharp},\cL)$, and it is 
an irreducible 
$\cG(P,\cL)$-module of weight one, that is by definition, 
any $\cG(P,\cL)$-submodule of $\Gamma(P,\cL)$ of weight
one is of the form 
$J\Gamma(P,\cL)$ for some ideal $J$ of $R$,
\item 
$\Gamma(P_0,\cL_0)$
is an irreducible 
$\cG(P_0,\cL_0)$-module of weight one. 
\end{enumerate}
\end{lemma}

See \cite[V, 2.4.2; VI, 1.4.7]{MB85} and \cite[Lemma~5.12]{Nakamura99}.

\begin{lemma}\label{lemma:K(P0,L0)} We define a morphism 
$\lambda(\cL_0):G^{\sharp}_0\to \op{Pic}^0(P_0)$ by 
$$\lambda(\cL_0)(a)=T_a^*(\cL_0)\otimes \cL_0^{-1}$$ 
for any $U$-valued point $a$ of $G^{\sharp}_0$, 
$U$ any $k(0)$-scheme.
Then 
\begin{enumerate}
\item $K(P_0,\cL_0)=\ker\lambda(\cL_0)$,
\item $\cG(P_0,\cL_0)$ is determined uniquely by $(P_0,\cL_0)$.
\end{enumerate} 
\end{lemma}
\begin{pf} 
 First we
note that $G^{\sharp}_0\simeq G_0\times (X/Y)$ in general,  
and that in the totally degenerate case $G_0\simeq 
\Hom_{\bZ}(X,\bG_m)$, while in the general case $G_0$ is a 
$\Hom_{\bZ}(X,\bG_m)$-torsor 
over an abelian variety $A_0$ whose extension class 
is determined uniquely by $(P_0,\cL_0)$. \par
Next we recall $\op{Pic}^0(P_0)$ is a $k(0)$-scheme 
by \cite[232, Corollaire 6.6]{FGA}. Then the proof 
of the first assertion is the same as 
\cite[Lemma~5.14]{Nakamura99}.  
\par
Next we prove the second assertion. 
We see as in the case of abelian varieties 
that $K(P_0,\cL_0)$ is the maximal 
subscheme of $G^{\sharp}_0$ such that the sheaf 
$m^*(\cL)\otimes p_2^*(\cL)^{-1}$ is trivial on $K(P_0,\cL_0)\times P_0$
where 
$m : G^{\sharp}_0\times P_0\to P_0$ is 
the action of $G^{\sharp}_0$ and 
$p_2:G^{\sharp}_0\times P_0\to P_0$ is the second projection. 
We also see as in the case of abelian varieties that  
$\cG(P_0,\cL_0)$ is the scheme representing 
the functor similar to 
\ref{subsec:Heisenberg group over k(eta)} :
\begin{equation*}
\begin{aligned}\label{eq:calG(L0)}
\cG(P_0,\cL_0)(U)&=\{(a,\phi);a\in K(P_0,\cL_0)(U)\ \op{and}\\
&\hskip 1cm \phi\ 
\op{is\ an\ isomorphism\ of}\  
T^*_a(\cL_{0,U})\ \op{with}\ \cL_{0,U}\}
\end{aligned}
\end{equation*}
for any $k(0)$-scheme $U$.
By the first assertion and \ref{eq:calG(L0)} 
$K(P_0,\cL_0)$ and $\cG(P_0,\cL_0)$ are 
independent of the choice of a Delaunay $g$-cell $\sigma$. 
\qed\end{pf}

\begin{defn} Let $k$ be an algebraically closed field and  
$(P_0,\cL_0)$ be a TSQAS over $k=k(0)$. 
Then we define 
\begin{align*}
e_{\min}(K(P_0,\cL_0))&=\max\{n>0; 
\ker(n\cdot\op{id}_{G^{\sharp}_0})\subset K(P_0,\cL_0)\},\\
e_{\max}(K(P_0,\cL_0))&=\min\{n>0; 
\ker(n\cdot\op{id}_{G^{\sharp}_0})\supset K(P_0,\cL_0)\}.
\end{align*}
where $G^{\sharp}_0$ is the closed fibre of $G^{\sharp}$ 
in \ref{subsec:group scheme G and Gsharp}.
If the order of $K(P_0,\cL_0)$ 
and the characteristic of $k(0)$ are coprime, 
then $K(P_0,\cL_0)\simeq\oplus_{i=1}^g(\bZ/e_i\bZ)^{\oplus 2}$ 
for some positive integers $e_i$ with $e_i|e_{i+1}$.  From this 
it follows 
$e_{\min}(K(P_0,\cL_0))=e_1$ and 
$e_{\max}(K(P_0,\cL_0))=e_g$.
\end{defn}

\begin{thm}\label{thm:psqas very ample}
Let $(Q_0,\cL_0)$ be 
a projectively stable quasi-abelian scheme over $k(0)$.  
If 
$e_{\min}(K(P_0,\cL_0))\geq 3$, 
$\Gamma(Q,\cL)\otimes k(0)$ is very ample on $Q_0$.
\end{thm}
\begin{pf} See \cite[Theorem~6.3]{Nakamura99}. 
\qed\end{pf}

\begin{cor}\label{cor:freeness of Gamma(P,L)} 
Suppose $e_{\min}(K(P_0,\cL_0))\geq 3$. 
Then $\Gamma(P,\cL)$ is base-point free and 
defines a finite morphism $\Phi_{\cL}$ of $P$ 
into the projective space.
The image of $P$ by $\Phi_{\cL}$ is isomorphic to $Q$. 
\end{cor}
\begin{pf} Let $\nu : P\to Q$ be the normalization morphism.
By definition $\cL_P=\nu^*(\cL_Q)$. 
By Lemma~\ref{lemma:comparison of Gamma} we have 
$\Gamma(P,\cL_P)=\nu^*\Gamma(Q,\cL_Q)$. Hence $\Gamma(P,\cL_P)$ 
is base-point free by Theorem~\ref{thm:psqas very ample} 
so that it defines a finite $S$-morphism
$\Phi_{\cL_P}:P\to\bP(\Gamma(P,\cL_P))$, 
which factors through $Q$. Since by Theorem~\ref{thm:psqas very ample} 
$\Gamma(Q,\cL)$ is very ample on $Q$,
the image $\Phi_{\cL_P}(P)$ is 
$Q$. 
\qed\end{pf}

%%%%%%%%%%%
%%%%%%%%%%%
\section{Level $G(K)$-structure}
\label{sec:level G(K) structure}
%%%%%%%%%%
%%%%%%%%%%
Let $\zeta_N$ be a primitive 
$N$-th root of unity and $\cO_N:=\bZ[\zeta_N, 1/N]$. 
For simplicity we write
$\cO=\cO_N$.

\begin{defn}\label{defn:G(K) and V(K)} 
Let $H$ be a finite abelian group such that  
$e_{\max}(H)$, the maximal order of elements in $H$, 
is equal to $N$.  
Now we regard $H$ as 
a constant finite abelian group 
$\cO$-scheme.
Let $H^{\vee}:=\Hom_{\cO}(H,\bG_{m,\cO})$ 
be the Cartier dual of $H$.  We set 
$K:=K(H)=H\oplus H^{\vee}$ and  
define a bimultiplicative (or simply a {\it bilinear}) form 
$e_K : K\times K\to \bG_{m,\cO}$ by 
\begin{equation*}
e_K(z\oplus\alpha,w\oplus\beta)=\beta(z)\alpha(w)^{-1}
\end{equation*} 
where $z,w\in H$, $\alpha,\beta\in H^{\vee}$.  
We note that $H$ 
is a maximally isotropic subgroup of $K$,  
unique up to isomorphism. 
\par
Let $\mu_{N}:=\op{Spec} \cO[x]/(x^{N}-1)$ be 
the group scheme of  
all $N$-th roots of unity. We define group $\cO$-schemes
$\cG(K)$ and $G(K)$ by 
\begin{align*}
\cG(K)&:
=\{(a,z,\alpha); a\in\bG_{m,\cO}, z\in H,\alpha\in H^{\vee}\},\\ 
G(K)&:
=\{(a,z,\alpha); a\in\mu_{N}, z\in H,\alpha\in H^{\vee}\}
\end{align*}
endowed with 
group scheme structure
\begin{equation*}
(a,z,\alpha)\cdot (b,w,\beta)
=(ab\beta(z),z+w,\alpha+\beta).
\end{equation*}

Let $V(K)$ be the group algebra $\cO[H^{\vee}]$ 
of $H^{\vee}$ over $\cO$, and 
an $\cO$-basis $v(\chi)$ $(\chi\in H^{\vee})$ of $V(K)$. 
The group scheme $G(K)$ and $\cG(K)$ 
act on $V(K)$ by 
\begin{equation*}U(K)(a,z,\alpha)(v(\chi))
=a\chi(z)v(\chi+\alpha)
\end{equation*}
where $a\in\mu_{N}$ or $a\in\bG_{m,\cO}$, 
$z\in H$ and  
$\alpha\in H^{\vee}$. Let ${\bar G}(K)=U(K)G(K)$.
\end{defn}
\begin{defn}\label{defn:weight one}
Let $k$ be a field over $\cO$. 
Any $G(K)\otimes k$-module $V$ is of weight one 
if every $a\in\mu_N\subset G(K)\otimes k$ 
acts on $V$ as scalar multiplication $a\cdot\op{id}_V$. 
Then we say that the action of $G(K)$ on $V$ is of
weight one. 
\end{defn}
\begin{rem}\label{rem:weight one, Mumford theorem}
By \cite{Mumford66}
$V(K)\otimes k$ is an irreducible $G(K)\otimes k$-module of weight one, 
unique up to equivalence. Any finite dimensional $G(K)$-module 
of weight one over $k$ is 
a direct sum of copies of $V(K)\otimes k$.\par
Let $R$ be a discrete valuation ring, $k(0)=R/I$ and $S=\Spec R$.
If the order of $K(P,\cL)$ and 
the characteristic of $k(0)$ are coprime, then 
$K(P,\cL)$ is a reduced flat finite group $S$-scheme, 
\'etale over $S$. Hence by taking a finite base change 
if necessary, we may assume 
by \cite[Section~7]{Nakamura99} that
$(K(P,\cL),e_S^{\sharp})\simeq (K_S,e_{K,S})$ and
$\cG(P,\cL)\simeq \cG(K)_S$ for a suitable $K$.\end{rem}
\begin{defn}\label{defn:finite Heisenberg group}
The (finite) Heisenberg group scheme $G(P,\cL)$ 
of $(P,\cL)$ is defined to be the unique subgroup scheme of 
$\cG(P,\cL)$ mapped isomorphically onto 
$G(K)_S$ when $\cG(P,\cL)\simeq\cG(K)_S$. 
See \cite[Section~7]{Nakamura99}.
Let $G(P_0,\cL_0)=G(P,\cL)\otimes k(0)$.
\end{defn}
\begin{defn}\label{defn:defn of K-symplectic tsqas}   
A pair $(P_0,\cL_0)$ is called 
a $g$-dimensional $K$-symplectic torically  
stable quasi-abelian scheme 
or a $K$-symplectic TSQAS over $k$ if
\begin{enumerate}
\item[(i)] $(P_0,\cL_0)$ is a $g$-dimensional torically 
stable quasi-abelian scheme over $k$, that is, a closed fibre of 
some $(P,\cL)$ in Theorem~\ref{thm:construction of P},
\item[(ii)] 
$(K(P_0,\cL_0),e^{\sharp}_{S,0})\otimes\bar k
\simeq (K,e_K)\otimes\bar k$. 
\end{enumerate}
\end{defn}

 Suppose $e_{\min}(K)\geq 3$ in what follows.
\begin{defn}\label{defn:level G(K) str} Let 
$(Z,L)$ a $K$-symplectic TSQAS over $k$. 
{\it A level $G(K)$-structure} $(\phi,\rho)$ 
on $(Z,L)$ is defined to be a pair of a finite $k$-morphism
$\phi:Z\to\bP(V(K)\otimes k)$ and 
a central extension isomorphism 
$\rho:G(K)\otimes_{\cO} k\to G(Z,L)$ 
which induce a weight one $\rho$-equivariant isomorphism of 
line bundles (namely the centers acting 
as scalar multiplication of weight one)
\begin{enumerate}
\item[(i)] $\phi^*(O_{\bP(V(K)\otimes k)}(1))\simeq L$.
\end{enumerate}
\par
\medskip
The isomorphism (i) induces
a weight one $\rho$-equivariant isomorphism of $k$-vector spaces 
by Lemma~\ref{lemma:comparison of Gamma} 
\begin{equation*}
H^0(\phi^*) : H^0(O_{\bP(V(K))}(1))\otimes k
=V(K)\otimes k\simeq H^0(Z,L)
\end{equation*}
where $G(K)$ acts on $V(K)\otimes_{\cO} k$ 
as a conjugate to $U(K)$. \par
\medskip
The level $G(K)$-structure $(\phi,\rho)$ 
is called {\it a rigid $G(K)$-structure} if 
\begin{enumerate}
\item[(ii)] $\rho=G(H^0(\phi^*))\circ (U(K)\otimes_{\cO} k)$
\end{enumerate}
where $G(H^0(\phi^*))(\bar g)
:=H^0(\phi^*)\circ\bar g\circ (H^0(\phi^*))^{-1}$ for any 
$\bar g\in {\bar G}(K)$, that is, 
the action of $G(K)$ on $V(K)\otimes_{\cO} k$ 
is just $U(K)$.\par
\medskip
If (i) is satisfied, we denote $(Z,L,G(Z,L),\phi,\rho)$ 
by $(Z,\phi,\rho)_{\LEV}$ because $L$ and $G(Z,L)$ 
are uniquely determined by $\phi$ and $\rho$. If (i) and 
(ii) are true, we denote it by 
$(Z,\phi,\rho)_{\RIG}$.
\end{defn}

\begin{defn}\label{defn:isomorphism of level structure}  
Let $(Z_i,L_i,G(Z_i,L_i),\phi_i,\rho_i)$ be $k$-TSQASes 
with level $G(K)$-structure $(i=1,2)$.  
They are 
defined to be isomorphic 
if there is a $k$-isomorphism $f : Z_1\simeq Z_2$ 
such that 
\begin{enumerate}
\item[(i)] $L_1\simeq f^*L_2$, 
\item[(ii)] $\rho_1=G(f^*)\circ
\rho_2$
\end{enumerate}
where $G(f^*)(g)=f^*g(f^*)^{-1}$ for any 
$g\in G(Z_2,L_2)$. \par
In this case we write 
$(Z_1,\phi_1,\rho_1)_{\LEV}\simeq (Z_2,\phi_2,\rho_2)_{\LEV}$. 
 By (ii) we have  
$G(Z_1,L_1)=G(f^*)G(Z_2,L_2)$. 
We define $(Z_i,\phi_i,\rho_i)_{\RIG}$ to be isomorphic 
if $(Z_i,\phi_i,\rho_i)_{\LEV}$ are isomorphic.
\end{defn}

\begin{lemma}\label{lemma:isom of rigid structure} 
Let  $(Z_i,\phi_i,\rho_i)_{\RIG}$ be $k$-TSQASes\ 
with rigid $G(K)$-structure $(i=1,2)$.    
Then the following are equivalent:
\begin{enumerate}
\item 
$(Z_1,\phi_1,\rho_1)_{\RIG}\simeq
(Z_2,\phi_2,\rho_2)_{\RIG}$,
\item there is a $k$-isomorphism
$f : Z_1\simeq Z_2$ with 
$\phi_1=\phi_2\circ f$.
\end{enumerate}
\end{lemma}
\begin{pf} Though the definition of rigid $G(K)$-structure is slightly
different from 
\cite[Section 9]{Nakamura99}, the proof of this lemma proceeds in the same
manner as 
\cite[Lemma 9.7]{Nakamura99}.\qed\end{pf}

\begin{lemma}\label{lemma:existence of rigid structure} 
Let $(Z,L)$ be 
a $K$-symplectic TSQAS over $k$. For 
any level $G(K)$-structure $(\phi,\rho)$ on $(Z,L)$ 
there exists 
a rigid $G(K)$-structure 
$(\phi(\rho),\rho)$ 
such that  
$(Z,\phi(\rho),\rho)_{\LEV}\simeq 
(Z,\phi,\rho)_{\LEV}$.
If $L$ is very ample, 
then the rigid $G(K)$-structure $(\phi(\rho),\rho)$
is unique.
\end{lemma}
\begin{pf}  We choose 
and fix a symplectic isomorphism 
$\sigma:(K,e_K)\otimes k(0)\simeq 
(K(P_0,\cL_0),e^{\sharp}_{S,0})$ and 
a central extension isomorphism 
$\rho:G(K)\otimes k(0)\simeq G(P_0,\cL_0)$ 
in a compatible way. We may assue $k=k(0)$.
By the uniqueness of weight one $G(K)\otimes k$-module
there is a $k$-isomorphism
$H^0(\phi(\rho)^*):V(K)\otimes k\to \Gamma(P_0,\cL_0)$ 
such that 
\begin{equation*}
\rho(g)(H^0(\phi(\rho)^*)(w))
=H^0(\phi(\rho)^*) U(K)(g)(w)\ (\forall g\in G(K), 
w\in V(K))
\end{equation*} 
Hence $\rho=G(H^0(\phi(\rho)^*))U(K)$.
By Corollary~\ref{cor:freeness of Gamma(P,L)}, 
$\phi(\rho)$ is a finite $k$-morphism 
of $Z$ into $\bP(V(K)\otimes_{\cO} k)$. 
Uniqueness of $\phi(\rho)$ follows from irreducibility of $U(K)$ 
and Schur's lemma when $L$ is very ample.
\qed\end{pf}

We note that if $L$ is not very ample, 
then there might be an automorphism $\psi$ of $(Z,L)$ 
which keeps $\phi$ and $\rho$ invariant.

\begin{defn}\label{defn:T-SQAS} Let $T$ be 
a noetherian $\cO$-scheme. Then a quintuplet\\
$(P,\cL,G(P,\cL),\phi,\rho)$ is 
called a torically stable quasi-abelian $T$-scheme 
(abbr. a $T$-TSQAS)
of relative dimension $g$ 
with level $G(K)$-structure  if 
\begin{enumerate}
\item[(i)] 
$P$ is a proper flat $T$-scheme
with a relatively ample invertible sheaf $\cL$,
\item[(ii)] $\phi$ is a finite $T$-morphism of 
$P$ into $\bP(V(K)\otimes_{\cO} O_T)$,
\item[(iii)] $G(P,\cL)$ is a finite flat reduced 
group $T$-scheme acting on $(P,\cL)$,
\item[(iv)] $\rho:G(K)_T\to G(P,\cL)$ 
is an isomorphism of group $T$-schemes,
\item[(v)]$\phi^*
(O_{\bP(V(K))}(1)\otimes_{\cO}M)\simeq\cL$ is a
$\rho$-equivariant isomorphism (of line bundles on $P$)
for some weight one action of 
$G(K)_T$ on $O_{\bP(V(K))_T}(1)$,
\item[(vi)]  for any geometric point $s$ of $T$, 
$(P_s,\phi_s,\rho_s)$ is a TSQAS 
of dimension $g$ over $k(s)$ with level $G(K)$-structure
\end{enumerate}
where $M$ is some invertible $O_T$-module with 
trivial $G(K)_T$-action and 
$\pi:P\to T$ is the structure
morphism.\par
\medskip

It follows from  the condition (vi)
that the isomorphism (v) induces a $\rho$-equivariant isomorphism 
of $O_T$-Modules
\begin{equation*}
R^0\pi_*(\phi^*)\ :\ V(K)\otimes_{\cO}M\ \simeq\ \pi_*(\cL).
\end{equation*} 
 
%\medskip

We denote $(P,\cL,G(P,\cL),\phi,\rho)$ 
by $(P,\phi,\rho)_{\op{LEV}}$. 
We call 
$(P,\phi,\rho)_{\op{LEV}}$ a $T$-TSQAS with rigid 
$G(K)$-structure and
denote it by $(P,\phi,\rho)_{\RIG}$ if
\begin{enumerate}
\item[(vii)] $R^0\pi_*(\rho)=G(R^0\pi_*(\phi^*))\circ U(K)_T$.
\end{enumerate} 
\end{defn}

\begin{defn}\label{defn:isom of T TSQAS} 
Let $(P_i,\phi_i,\rho_i)_{\op{LEV}}:=(P_i,\cL_i,G(P_i,\cL_i),\phi_i,\rho_i)$ 
 $(i=1,2)$ 
be $T$-TSQASes 
with level $G(K)$-structure.  We define 
them to be isomorphic
if there exist a $T$-isomorphism 
$f:P_1\to P_2$ and an invertible $O_T$-module 
$M$ with 
trivial $G(K)_T$-action
such that 
\begin{enumerate}
\item[(i)] $\cL_1\simeq f^*\cL_2\otimes_{\cO}M$,
\item[(ii)] 
$\rho_1=G(f^*)\circ\rho_2$.
\end{enumerate}

If $(P_i,\phi_i,\rho_i)$ are rigid $G(K)$-structures, then
(i) and (ii) is by Lemma~\ref{lemma:isom of rigid structure} 
equivalent to 
\begin{enumerate}
\item[(iii)] $\phi_1=\phi_2\circ f$.
\end{enumerate}

\medskip
An algebraic space $T$ is by definition 
the isomorphism class of  
an \'etale representative $U\to T$ 
with \'etale equivalence relation $R\subset U\times U$. See \cite{K71}. 
Let $p_i: R\to U$ be the composite of 
the immersion $R\subset U\times U$ with $i$-th projection $(i=1,2)$. 
A $T$-TSQAS $(Z, \psi, \rho)_{\LEV}$ with level $G(K)$-structure
is  a $U$-TSQAS $(Z_U, \psi_U, \rho_U)_{\LEV}$  
whose pullbacks by $p_i$ are isomorphic as $R$-TSQASes with level
$G(K)$-structure. 
\end{defn}

\begin{defn}
We define the functor $\cal{SQ}^{toric}_{g,K}$ as follows. 
For any noetherian ${\cO}_N$-scheme $T$, we set
\begin{align*}
\cal{SQ}^{toric}_{g,K}(T)&
=\op{the\ set\ of\ torically\ stable\ 
quasi-abelian}\\ 
&\hskip 0.5cm \ T\op{-schemes}\ (P,\phi,\rho)_{\op{LEV}}
\ \op{of\ relative\ dimension}\, g\\
&\hskip 0.5cm  
\, \op{with\ level}\,G(K)\op{-structure}\ 
\op{modulo}\ T\op{{\rm -}isom}.
\end{align*}
\end{defn}
\subsection{General rigid $G(K)$-structures}
\label{susbsec:General rigid G(K) structures}
Even if $(Z,L)$ is neither a TSQAS nor a PSQAS, we can also speak of 
a rigid (or level) $G(K)$-structure. 
If $G(K)$ acts on a polarized $k$-scheme $(Z,L)$ 
with $L$ very ample (or equivalently, $L$ is $G(K)$-linearized 
\cite[p.~30]{MFK94}),
then $H^0(Z, L)$ becomes a $G(K)$-module in a natural manner. Let $\rho$
be the action of $G(K)$ on $H^0(Z, L)$ and $\bP:=\bP(V(K)\otimes_{\cO} k)$.
If $\phi: (Z,L)\to (\bP, O_{\bP}(1))$ is 
a $G(K)$-equivariant closed immersion such that $H^0(\phi^*) :
V(K)\otimes_{\cO} k\to H^0(Z,L)$ is an
isomorphism with $\rho=G(H^0(\phi^*))\circ (U(K)\otimes k)$, 
then we call the triplet
$(Z,\phi,\rho)$ a rigid $G(K)$-structure on $(Z,L)$, 
which we denote $(Z,\phi,\rho)_{\RIG}$. 
In particular, if $Z$ is a $G(K)$-invariant closed subscheme of $\bP$, then 
$(Z,i,U(K))$ is a rigid $G(K)$-structure on $(Z,O_{\bP}(1)\otimes O_Z)$ where 
$i$ is the inclusion of $Z$ in $\bP$. If $\rho$ is of weight one, 
then Lemmas~\ref{lemma:isom of rigid structure},
\ref{lemma:existence of rigid structure} are true as well. 
%%%%%%%%%%
%%%%%%%%%%
\section{The stable reduction theorem}
\label{sec:The stable reduction theorem}
%%%%%%%%%%%%%%%%%
%%%%%%%%%%%%%%%%%
\subsection{The rigid $G(K)$-structure we start from}
\label{subsec:rigid str we start}
Let $R$ be a complete 
discrete valuation ring, 
$k(\eta)$ (resp. $k(0)$) the fraction field  
(resp.  the residue field) of $R$, and $S=\Spec R$.  Let 
$(G_{\eta}, \cL_{\eta})$ be
a polarized abelian variety over $k(\eta)$ 
with $\cL_{\eta}$ ample 
and $K(\cL_{\eta}):=\ker\lambda(\cL_{\eta})$.
Let $e^{\cL_{\eta}}$ be the Weil pairing of 
$K(\cL_{\eta})$. 
\par 
Suppose that the order of $K(\cL_{\eta})$ and 
the characteristic of $k(0)$ are coprime.
Then there exists a finite symplectic constant 
abelian group $\bZ$-scheme
$(K,e_K)$ such that 
$(K,e_K)\otimes_{\bZ} k(\eta)
\simeq (K(\cL_{\eta}),e^{\cL_{\eta}})$.
Let $N=e_{\max}(K)$.
We also may assume  
that $R$ contains 
a primitive $N$-th root $\zeta_N$ of unity. If  
$e_{\min}(K)\geq 3$, then  
by Lemma~\ref{lemma:existence of rigid structure} 
$(G_{\eta},\cL_{\eta})$ 
has a unique rigid $G(K)$-structure   
because $\cL_{\eta}$ is very ample.

\begin{thm}\label{thm:summary}
Let $R$ be a complete discrete valuation ring and $S=\Spec R$.  
Let $(G_{\eta}, \cL_{\eta})$ be
a polarized abelian variety over $k(\eta)$. Let 
$K(\cL_{\eta}):=\ker\lambda(\cL_{\eta})$ and 
$N=e_{\max}(K(G_{\eta},\cL_{\eta}))$. 
Assume that 
\begin{enumerate}
\item[(i)] the characteristic of $k(0)$ and 
the order of $K(\cL_{\eta})$ are coprime, 
\item[(ii)] $e_{\min}(K(\cL_{\eta}))\geq 3$,
\item[(iii)] 
$R$ contains a primitive $N$-th root $\zeta_N$ of unity. 
\end{enumerate}

Then after a suitable finite base change 
if necessary, there exist 
flat projective schemes 
$(P, \cL)$ and $(Q,\cL)$, 
semiabelian group schemes $G$ and $G^{\sharp}$, 
the flat closure $K(P,\cL)$ 
of $K(\cL_{\eta})$ in $G^{\sharp}$, 
a symplectic form 
$e^{\sharp}_S$ on $K(P,\cL)$ extending 
$e^{\cL_{\eta}}$
and the Heisenberg group 
schemes $G(P,\cL)$ and $\cG(P,\cL)$ of $(P,\cL)$, 
all of these being defined over $S$,
such that 
\begin{enumerate}
\item\label{item:reduced} $P$ is reduced over $S$,
\item\label{item:open subgroup} $(G,\cL)$ 
and $(G^{\sharp},\cL)$ are 
open subschemes of 
$(P,\cL)$,
\item\label{item:Gsharp} $G^{\sharp}=K(P,\cL)\cdot G$,
\item\label{item:generic fibre} $(G_{\eta},\cL_{\eta})
\simeq (G^{\sharp}_{\eta},\cL_{\eta})
\simeq (P_{\eta},\cL_{\eta})\simeq 
(Q_{\eta},\cL_{\eta})$,
\item\label{item:isomorphisms} there exists a constant finite
symplectic  abelian group $\bZ$-scheme $(K,e_K)$ such that 
$(K(P,\cL),e^{\sharp}_S)\simeq (K,e_K)_S$ 
and $G(P,\cL)\simeq G(K)_S$,  
\item\label{item:irreducible representation} 
$\Gamma(G^{\sharp},\cL)\simeq\Gamma(P,\cL)
\simeq\Gamma(Q,\cL)\simeq V(K)\otimes_{\cO_N}R$ and they are 
irreducible $G(P,\cL)$-modules of weight one, 
unique up to equivalence.
\end{enumerate}
\end{thm}
 
See \cite{Nakamura99} 
for the proof of it and for the details of $(Q,\cL)$.

%%%%%%%%%%%%%%%%%%%%%%
%%%%%%%%%%%%%%%%%%%%%%
%%%%%%%%%%%%%%%%%%%%%%
\section{The scheme parametrizing TSQASes}
\label{sec:The scheme parametrizing TSQASes}
%%%%%%%%%%%%%%%%%%%%%%
%%%%%%%%%%%%%%%%%%%%%%
Let $K$ be a symplectic finite abelian group with 
$e_{\min}(K)\geq 3$. Let $N=e_{\max}(K)$ and $\cO=\cO_N$. 
In what follows we fix $K$ and $\cO$.

\subsection{The scheme $H_1\times H_2$}
\label{subsec:The scheme H1xH2}
Choose and fix a coprime pair of natural 
integers $d_1$ and $d_2$ such that
$d_1>d_2\geq 2g+1$ and $d_k\equiv 1\mod N$.  
This pair does exist because it is enough to
choose prime numbers $d_1$ and $d_2$ large enough 
such that $d_k\equiv 1\mod N$ and $d_1>d_2$. 
We choose integers $q_k$ such that $q_1d_1+q_2d_2=1$. 
\par
Now consider a $G(K)$-module $W_k(K):=W_k\otimes V(K)\simeq 
V(K)^{\oplus N_k}$ 
where $N_k=d_k^g$ and 
$W_k$ is a free $\cO$-module of rank $N_k$ 
with trivial $G(K)$-action.
Let $\rho_k$ be the natural action of $G(K)$ on $W_k(K)$.
\par
Let $H_k$ be the Hilbert scheme parametrizing all
pure $g$-dimensional polarized subschemes $(Z_k,M_k)$ 
of the projective space $\bP(W_k(K))$ such that
\begin{equation*} 
\chi(Z_k,nM_k)=n^gd_k^g\sqrt{|K|}\quad \op{for}\ k=1,2.
\end{equation*}

Let $X_k$ be the universal subscheme of $\bP(W_k(K))$ 
over $H_k$, $X$ the product of $X_1$ and $X_2$ over $\cO$. 
Hence $X$ is a subscheme of $\bP(W_1(K))\times_{\cO}\bP(W_2(K))$ 
over $H_1\times_{\cO} H_2$.

%%%%%%%%%%%%%%
%%%%%%%%%%%%%%
\subsection{The scheme $U_1$}
\label{subsec:The scheme U1}Let $(X,L)$ be a polarized 
$\cO$-scheme with $L$ very ample and  $P(n)$ an 
arbitrary polynomial.
Let $\Hilb^P(X)$ be the Hilbert scheme 
parametrizing all subschemes $Z$ of $X$ 
with $\chi(Z,nL_Z)=P(n)$. As is well known $\Hilb^P(X)$ is a projective 
$\cO$-scheme. \par
For a given projective scheme $T$ and 
a given flat projective $T$-scheme $(X,L)$ and an 
arbitrary polynomial $P(n)$, let 
$\Hilb^{P}_{\op{conn}}(X/T)$ be the scheme 
parametrizing all connected subschemes $Z$ of $X$ 
with $\chi(Z,nL_Z)=P(n)$
projected to one point of $T$. 
Then $\Hilb^{P}_{\op{conn}}(X/T)$
is a projective $T$-subscheme of $\Hilb^P(X)$.

Let $p_k:X_1\times_{\cO} X_2\to X_k$ 
be the $k$-th projection. Let $X=X_1\times_{\cO} X_2$ and 
$H=H_1\times_{\cO} H_2$. The aim of the subsequent sections 
is to construct a new compactification 
of the  moduli scheme of abelian varieties as the quotient 
of a certain subscheme of $\Hilb^{P}_{\op{conn}}(X/H)$ by 
$\GL(W_1)\times\GL(W_2)$.\par
Let $B$ be the pullback to $X$ 
of a very ample invertible sheaf on $H$.
Let $M_i=p_i^*(O_{\bP(W_i(K))}(1))$ and $M=d_2M_1+d_1M_2+B$.
Then $M$ is a very ample invertible sheaf on $X$.
Now we define $U_1$ to be the subset of 
$\Hilb^{P}_{\op{conn}}(X/H)$ consisting of all subschemes $Z$ such that 
\begin{enumerate}
\item[(i)] ${p_k}_{|Z}$ is an isomorphism  $(i=1,2)$,
\item[(ii)] $d_2L_1=d_1L_2$,
\item[(iii)] $Z$ is $G(K)$-stable
\end{enumerate}
where $P(n)=(2nd_1d_2)^g\sqrt{|K|}$ and $L_i=M_i\otimes O_Z$.\par
We prove that $U_1$ is a nonempty closed $\cO$-subscheme of 
$\Hilb^{P}_{\op{conn}}(X/H)$.\par
The condition $d_2L_1=d_1L_2$ is closed, while 
the condition that ${p_i}_{|Z}$ 
is an isomorphism is open and closed. The $G(K)$-stability of $Z$ 
is equivalent to the condition that 
$Z\in\Hilb^{P}_{\op{conn}}(X/H)$ is fixed by the natural $G(K)$-action 
on $X$ and $H$.
Hence it is a closed condition.
Hence $U_1$ is a closed, hence a projective $\cO$-subscheme of 
$\Hilb^{P}_{\op{conn}}(X/H)$.\par
It remains to show $U_1\neq\emptyset$. Let $k$ be an algebraically 
closed field over $\cO$, and $(A,L)$ 
a polarized abelian variety over $k$ 
with $G(A,L)\simeq G(K)$. Since $e_{\min}(K)\geq 3$, $L$ is 
very ample and $(A,d_iL)\in H_i$. 
Moreover $\Gamma(A,L)\simeq V(K)\otimes k$. Hence there is a unique 
rigid $G(K)$-structure of $(A,L)$ 
by Lemma~\ref{lemma:existence of rigid structure}, in other words, 
there is a unique $G(A,L)$-$G(K)$ equivariant closed immersion 
$\phi : A\to \bP(V(K)\otimes k)$ of $(A,L)$.
In particular $\phi(A)\ (\simeq A)$ is 
a $G(K)$-stable subscheme of $\bP(V(K)\otimes k)$. 
Since $L$ has a $G(A,L)$-linearization of weight one, 
$d_iL$ has a $G(A,L)$-linearization of weight $d_i$ too.
Since $d_i\equiv 1\mod N$, and since $a^N=1$ 
for any $a\in \mu_N$,
$d_iL$ has a $G(A,L)$-linearization of weight one.   
Hence $\Gamma(A, d_iL)$ is a direct sum of $V(K)\otimes k$.
Since $\Gamma(A,d_iL)$ is very ample, we can choose a 
$G(A,L)$-$G(K)$-equivariant closed immersion 
$\phi_i : A\to \bP(W_i(K))$. Then $\phi_i(A)$ is a $G(K)$-stable subscheme.
Therefore $(\phi_1(A),\phi_2(A))\in H_1\times H_2$. 
Let $Z\subset\phi_1(A)\times\phi_2(A)\simeq A\times A $ be the 
inverse image of the diagonal. 
Since $Z\simeq A$, we see that 
\begin{align*}
\chi(Z,&\,n(d_2L_1+d_1L_2+B)_Z)\\
&=\chi(A,2nd_1d_2L)=(2nd_1d_2)^g\sqrt{|K|}=P(n),
\end{align*}
which proves $Z\in\Hilb^{P}_{\op{conn}}(X/H)$. Thus $Z\in U_1(k)$. 
Hence $U_1\neq\emptyset$.
\par

\begin{lemma}\label{lemma:defn of L}
Let $k$ be an algebraically 
closed field over $\cO$. Let $Z\in U_1(k)$ and 
$L=q_1L_1+q_2L_2$. 
Then $L_i=d_iL$.
\end{lemma}
\begin{pf}One sees readily
$d_iL=d_i(q_1L_1+q_2L_2)=(d_1q_1+d_2q_2)L_i=L_i$.
\qed\end{pf}
%%%%%%%%%%%%%%
%%%%%%%%%%%%%%
\subsection{The scheme $U_2$}
\label{subsec:The scheme U2}
Let $q_i$ be the integers chosen above. 
Let $U_2$ be the open subscheme of $U_1$ 
consisting of all subschemes $Z$ of $X$ such that 
(i)-(iii) are true and 
\begin{enumerate}
\item[(iv)] $L$ is ample where $L=q_1L_1+q_2L_2$,
\item[(v)] $\chi(Z,nL)=n^g\sqrt{|K|}$, 
\item[(vi)] $H^q(Z,nL)=0$ 
for $q>0$ and $n>0$,
\item[(vii)] $\Gamma(Z,L)$ is base point free,
\item[(viii)] the pullback  
$H^0(p_i^*) : W_i(K)\otimes k\to\Gamma(Z,L_i)$ 
is surjective for $i=1,2$,
\item[(ix)] $Z$ is reduced.
\end{enumerate}

It is clear that (iv)-(ix) are 
open conditions. Note that the
surjectivity in (viii) implies a $G(K)$-equivariant isomorphism 
in view of (iii).\par
We note $U_2\neq\emptyset$. In fact, letting $k$ be 
an algebraically closed field over $\cO$ we
choose a polarized abelian variety $(A,L)$ over $k$
with $G(A,L)\simeq G(K)$. Then $L$ is 
very ample and $(A,d_iL)\in H_i$ (identified with $\phi_i(A)$), and
the inverse image $Z$ of the diagonal $(\simeq A)$ 
belongs to $U_1(k)$ as we saw 
in \ref{subsec:The scheme U1}. 
Since $L_i=d_iL$ by 
Lemma~\ref{lemma:defn of L}, 
all the conditions (iv)-(viii) are true  
for $Z$ as is well known. Hence $Z\in U_2(k)$. Hence $U_2\neq\emptyset$.
\subsection{The scheme $U_3$}\label{subsec:The scheme V3}
Next we recall that the locus $U_{g,K}$ 
of abelian varieties is an open subscheme of $U_2$. 
The condition on $Z\in U_{g,K}$ 
that the natural action of $G(K)$ on $Z$ 
is contained in $\Aut^0(Z)$  is open. 
Since $L$ is very ample, 
this condition implies that if $(Z,L)$ is a polarized abelian variety, then
the restriction of the $G(K)$-action to $Z$ reduces to $K(Z,L)$.\par
Now we define $U_3$ to be a reduced subscheme in $U_2$ 
whose underlying set is 
the union of all the irreducible components of $U_2$ 
over which at least one of the  
geometric fibres of $X$ is a polarized abelian variety 
$(Z,L)$ with $L_i=d_iL$ and $L$ very ample such that
\begin{enumerate}
\item[(x)]  the restriction to $Z$ of the $G(K)$-action on $X$ 
is contained in $\Aut^0(Z)$.
\end{enumerate} 

By definition $U_3$ is the closure of $U_{g,K}$ in $U_2$ 
with reduced structure. Note that it is  
an $\cO$-subscheme of $\Hilb^{P}_{\op{conn}}(X/H)$. 

%%%%%%%%%%%%%%
%%%%%%%%%%%%%%

\section{The fibres over $U_3$}
\label{sec:The fibres over U3}
\subsection{(S$_k$) and (R$_k$)}\label{subsec:Sk and Rk}
Here we recall the conditions (S$_k$) and (R$_k$): 
\begin{align*}
(\op{S}_k) &\qquad \depth(A_p)\geq \inf(k,\op{ht}(p))
\quad \op{for\ all}\ p\in\Spec(A),\\
(\op{R}_k)&\qquad A_p\ \op{is\ regular}\
\op{for\ all}\ p\in\Spec(A)\ \op{with}\ \op{ht}(p)\leq k\ .
\end{align*}

\begin{lemma}\label{lemma:normality reducedness} 
Let $A$ be a noetherian local ring. Then
\begin{enumerate}
\item $\op{(Serre)}$\ $A$ is normal if and only if 
$(\op{R}_1)$ and $(\op{S}_2)$ are true for $A$,
\item $A$ is reduced if and only if 
$(\op{R}_0)$ and $(\op{S}_1)$ are true for $A$.
\end{enumerate}
\end{lemma} 
See \cite[Theorem~39]{Matsumura70} and   
\cite[IV$_2$, 5.8.5 and 5.8.6]{EGA}.

\begin{lemma}\label{lemma:valuative lemma for separatedness}
Let $R$ be a discrete valuation ring, $S:=\Spec R$,
$\eta$ the generic point of $S$ and 
$k(\eta)$ the fraction field of $R$. Let $(Z_k,\phi_k,\rho_k)_{_{\RIG}}$ 
be flat proper schemes over $S$
with rigid $G(K)$-structure in the sense of 
\ref{susbsec:General rigid G(K) structures} such that 
$\phi_k$ is a closed immersion of $Z_k$ into $\bP(V(K)\otimes_{\cO}R)$ and 
$H^0(\phi_k^*): V(K)\otimes_{\cO}R\to H^0(Z_k,L_k)$ is a $G(K)$-isomorphism.
If $(Z_k,\phi_k,\rho_k)_{_{\RIG}}$ are 
isomorphic abelian varieties over $k(\eta)$, 
then they are isomorphic over $S$.
\end{lemma}
\begin{pf} 
Let $H=\Hilb^{P(n)}_{\bP(V(K))}$ be the Hilbert scheme parametrizing all 
subschemes of $\bP(V(K))$ with Hilbert polynomial $P(n)=n^g\sqrt{|K|}$ and 
$X_{\op{univ}}$ the universal subscheme of $\bP(V(K))$ over $H$.
Then $\phi_k$ induces a unique morphism $\Hilb(\phi_k):S\to H$ such that 
$Z_k$ is the pullback by $\Hilb(\phi_k)$ of $X_{\op{univ}}$. 
By the assumption there is a $k(\eta)$-isomorphism 
$f_{\eta} : Z_{1,\eta}\to Z_{2,\eta}$ by Lemma~\ref{lemma:isom of rigid
structure}
such that $\phi_{1,\eta}=\phi_{2,\eta}\circ f_{\eta}$.
It follows that $\Hilb(\phi_{1,\eta})=\Hilb(\phi_{2,\eta})$. 
Since $H$ is separated, $\Hilb(\phi_1)=\Hilb(\phi_2)$, 
hence $\phi_1(Z_1)=\phi_2(Z_2)$. 
This implies that there is an $S$-isomorphism $f:Z_1\to Z_2$ 
extending $f_{\eta}$ such that $\phi_1=\phi_2\circ f$. 
This proves 
$(Z_1,\phi_1,\rho_1)_{_{\RIG}}\simeq 
(Z_2,\phi_2,\rho_2)_{_{\RIG}}$ by Lemma~\ref{lemma:isom of rigid
structure}.
\qed\end{pf}

\begin{lemma}\label{lemma:geom fibres over U3}
Let $R$ be a discrete valuation ring, $S:=\Spec R$,
$\eta$ the generic point of $S$ and 
$k(\eta)$ the fraction field of $R$. Let 
$h$ be a morphism from $S$ into $U_3$.
Let $(Z,\cL)$ be the pullback by $h$ of the universal subscheme 
$Z_{\op{univ}}$, universal for $\Hilb^{P}_{\op{conn}}(X/H)$, such that 
$(Z_{\eta},\cL_{\eta})$ 
is a polarized abelian variety. 
Then $(Z,\cL)$ is isomorphic to a (modified) Mumford's family 
$(P,\cN)$ in Theorem~\ref{thm:construction of P} 
after a finite base change if necessary.
\end{lemma}
\begin{pf} By the assumption on $h$, $(Z_{\eta},\cL_{\eta})$ 
is a polarized abelian variety over $k(\eta)$
such that the action of $G(K)$ is contained 
in $\Aut^0(Z_{\eta})\cap\Aut(Z_{\eta},\cL_{\eta})$. 
After a suitable finite base change 
we may assume by  
Theorem~\ref{thm:summary} 
that there is a flat projective family 
$(P,\cN)$ associated with the degeneration data 
of $(Z_{\eta},\cL_{\eta})$ such that 
\begin{enumerate}
\item $(P,\cN)$ is a $K$-symplectic $S$-TSQAS, 
in particular, $P_0$ is connected 
reduced and $P$ is normal,
\item $(P_{\eta},\cN_{\eta})\simeq (Z_{\eta},\cL_{\eta})$,
\item $G(P,\cN)\simeq G(K)_S$,
\item there is a $G(P,\cN)$-$G(K)_S$-equivariant polarized
finite morphism 
\begin{equation*}
\psi : (P, \cN)\to (\bP(V(K))_S, O_{\bP(V(K))_S}(1))
\end{equation*} 
such that $\psi_{\eta}$ is a closed immersion.
\end{enumerate}

Let $(Q,\cN_Q)$ be an $S$-PSQAS extending $(P_{\eta},\cN_{\eta})$ to $S$. 
The scheme $Q$ was defined in Section~\ref{sec:degenerating families}.
By Lemma~\ref{cor:freeness of Gamma(P,L)} 
$Q$ is the image of $P$ by the morphism 
$\psi : P\to \bP(V(K)\otimes_{\cO} R)$ 
defined by $\Gamma(P,\cN)$, while $\cN_Q$ is the restriction of 
$\cO_{\bP(V(K))}(1)_S$ to $Q$. Then
$\psi : P\to Q$ is the normalization of $Q$ 
in view of Theorem~\ref{thm:construction of P}.\par
Let $\pi : Z\to S$ be the flat family given at the start. 
Then $\Gamma(Z,\cL)$ is 
a free $R$-module of rank $\sqrt{|K|}$ by (ii). It is 
a $G(K)$-module of weight one, hence $G(K)$-isomorphic to $V(K)\otimes R$ 
after a finite base change.
See Remark~\ref{rem:weight one, Mumford theorem}.
By (vii) $\Gamma(Z,\cL)$ is base point free, which defines a 
$G(K)$-equivariant finite morphism $q : Z\to \bP(V(K)\otimes R)$ such that 
$q_{\eta}$ is a closed immersion 
of $Z_{\eta}$ because $e_{\min}(K)\geq 3$. 
Let $W$ be the flat closure of $q_{\eta}(Z_{\eta})$ in $\bP(V(K)\otimes
R)$, 
and $\cL_W$ the restriction of $\cO_{\bP(V(K))_S}(1)$. 
Since $Z_{\eta}$ is reduced, $W_{\red}$ is the flat closure 
of $q_{\eta}(Z_{\eta})$.  Hence $W$ is reduced.  
Since $Z_{\eta}$ is irreducible, so is $W$. 
Since $Z_0$ is reduced, so is $Z$, hence $q$ factors through $W$. 
It follows that $q:Z\to W$ is a finite surjective birational morphism.\par 
By \ref{susbsec:General rigid G(K) structures} and 
Lemma~\ref{lemma:existence of rigid structure} 
the $S$-scheme $(W,\cL_W)$ has a  
unique rigid $G(K)$-structure 
$(W,i_W,U(K))_{\RIG}$, while $(Q,\cN_Q)$ has a unique 
rigid $G(K)$-structure 
$(Q,i_Q,U(K))_{\RIG}$ by Theorem~\ref{thm:summary} 
where $i_W$ and $i_Q$ are natural inclusions of $W$ and $Q$ into 
$\bP(V(K))_S$.
Since we have
\begin{equation*}
(W_{\eta},i_{W_{\eta}},U(K))_{\RIG}\simeq  
(Q_{\eta},i_{Q_{\eta}},U(K))_{\RIG}
\simeq (Z_{\eta},i_{Z_{\eta}},\rho_{Z_{\eta}})_{\RIG},
\end{equation*}
$(W,i_W,U(K))_{\RIG}$ and $(Q,i_Q,U(K))_{\RIG}$ 
are $S$-isomorphic by Lemma~\ref{lemma:valuative lemma for separatedness}. 
It follows the action of $G(K)$ on $(W, \cL_W)$ is the same as 
that of $G(W, \cL_W)$.
\par
Hence to complete our proof of the theorem
it suffices to prove that $Z$ is also the normalization of 
$W$. Since $q$ is finite and birational, 
it suffices to prove that $Z$ is normal. 
By Lemma~\ref{lemma:normality reducedness} 
it suffices to check that (R$_1$) and (S$_2$) are true for $O_Z$. 
Since $Z_0$ is reduced, it is smooth at a generic point 
of any irreducible component of it. Hence $Z$ is smooth 
at any codimension one point of $Z$ supported by $Z_0$. 
Since $Z_{\eta}$ is smooth, $Z$ is 
codimension one nonsingular everywhere. 
This is (R$_1$). \par
Next we prove (S$_2$). Since $\pi:Z\to S$ is flat, 
any regular parameter $t$ of $R$ is not a zero divisor of $O_Z$. 
Let $p$ be a prime ideal of $O_Z$. 
If $p\cap R\neq 0$, then $t\in p$ 
and $q:=p/tO_Z$ is a prime ideal of $O_{Z_0}$ 
with $\op{ht}(q)=\op{ht}(p)-1$.
Since $Z_0$ is reduced, hence (S$_1$) for $Z_0$ is true by 
Lemma~\ref{lemma:normality reducedness}. It follows that 
 $\depth(O_Z)_p=\depth(O_{Z_0})_q+1
\geq\inf(1,\op{ht}(q))+1=\inf(2,\op{ht}(p))$. 
If $p\cap R=0$, then $k(\eta)\subset (O_Z)_p$ 
and $(O_Z)_p=(O_{Z_{\eta}})_{pO_{Z_{\eta}}}$. 
Hence $\depth(O_Z)_p=\dim(O_Z)_p\geq\inf(2,\op{ht}(p))$ 
because $Z_{\eta}$ is nonsingular. 
This proves (S$_2$).\par
Therefore $Z$ is normal by Lemma~\ref{lemma:normality reducedness}. 
It follows that $Z$ is the normalization of $W$ and 
$(Z,\cL)\simeq (P,\cN)$.
\qed\end{pf}

\begin{cor}\label{cor:closed fibre over U3} 
Let $(Z_0,\cL_0)$ be the closed fibre of $(Z,\cL)$ 
in Lemma~\ref{lemma:geom fibres over U3}. 
Then $(Z_0,\cL_0)$ is a $K$-symplectic TSQAS 
such that the action of $G(K)$ on $(Z_0,\cL_0)$ is $G(Z_0,\cL_0)$.
\end{cor}
\begin{pf} By the proof of 
Lemma~\ref{lemma:geom fibres over U3}, 
we see that $Z$ is the normalization of $W$ and 
$(W,i_W,\rho_W)_{\RIG}\simeq (Q,i_Q,\rho_Q)_{\RIG}$.  
The normalization morphism of $Z$ onto $W$ 
is $G(K)$-equivariant and
the action of $G(K)$ on $(W,\cL_W)$ is $G(W,\cL_W)$ by the proof of 
Lemma~\ref{lemma:geom fibres over U3}. Hence 
the action of $G(K)$ on $(Z,\cL)$ is $G(Z,\cL)$. This proves the corollary.
\qed\end{pf}

\begin{cor}\label{cor:fibres over U3 and G(Z,L)} 
Let $\Spec k$ be a geometric point over $\cO$ and 
$Z\in U_3(k)$. Let $L=M\otimes O_Z$ 
under the notation of \ref{subsec:The scheme U1}. 
Then $(Z,L)$ is a\/ $K$-symplectic TSQAS such that the $G(K)$-action on 
$(Z,L)$ induced from that on $W_i(K)$ is exactly $G(Z,L)$.
\end{cor}
\begin{pf} It follows from 
Corollary~\ref{cor:closed fibre over U3} that $(Z,L)$ is a\/ $K$-symplectic
TSQAS. 
The remaining assertion follows from Definition~\ref{defn:cal G(P,L)} 
and Theorem~\ref{thm:summary}.
\qed\end{pf}
%%%%%%%%%%%%%
%%%%%%%%%%%%%
\section{The geometric quotient}
\label{sec:The geometric quotient}
Lt $N=e_{\max}(K)$ and 
$\cO=\cO_N=\bZ[\zeta_N, 1/N]$. 
\begin{lemma}\label{lemma:GL(Wk)orbit}
Let $k$ be an algebraically closed field 
over $\cO$.
\begin{enumerate}
\item $U_3$ is $\GL(W_1)\times\GL(W_2)$-invariant, 
\item Let $(Z, L)\in U_3(k)$ and $(Z', L')\in U_3(k)$ where 
$L=M\otimes O_Z$ and $L'=M\otimes O_{Z'}$.  
If $(Z, L)\simeq (Z', L')$ 
as polarized varieties with $G(K)$-linearization, 
then $(Z', L')$ belongs to 
the $\GL(W_1)\times\GL(W_2)$-orbit of $(Z, L)$.
\end{enumerate}
\end{lemma}
\begin{pf}First we prove (2). 
 Let $f:(Z, L)\to (Z', L')$  be an isomorphism with $G(K)$-linearization. 
Hence $(Z, d_iL)$ and $(Z', d_iL')$ are isomorphic as 
polarized schemes with $G(K)$-linearization.
By the assumptions on $(Z, L)$ and $(Z', L')$, we see first 
$d_iL$ and $d_iL'$ are very ample. Hence 
$(Z, d_iL)$ and $(Z', d_iL')\in H_i(k)$ $(i=1,2)$.
Thus we see
\begin{enumerate}
\item  
there are commutative diagrams of $G(K)$-equivariant isomorphisms 
\begin{equation*}
\CD 
(Z, d_iL)@>{f_i}>> (Z', d_iL')\\
@VV{\iota_i}V @VV{\iota'_i}V\\ 
(\bP(W_i(K)), O_{\bP(W_i(K))}(1)) @>{F_i}>>
(\bP(W_i(K)), O_{\bP(W_i(K))}(1)).\\
\endCD
\end{equation*}
where $\iota_i$ and $\iota'_i$  are closed immersion of 
$(Z, d_iL)$ into $\bP(W_i(K))$.
\item there are commutative diagrams of $G(K)$-equivariant isomorphisms 
\begin{equation*}
\CD 
H^0(Z, d_iL)@<{H^0(f_i^*)}<< H^0(Z', d_iL')\\
@AA{H^0(\iota_i^*)}A @AA{H^0((\iota'_i)^*)}A\\ 
H^0(O_{\bP(W_i(K))}(1))\otimes k @<{H^0(F_i^*)}<<
H^0(O_{\bP(W_i(K))}(1)))\otimes k\\
\endCD
\end{equation*}
where $H^0(O_{\bP(W_i(K))}(1))=W_i(K)$. 
\end{enumerate}

Let $\rho_i$ and $\rho'_i$ be the $G(K)$-actions 
on $W_i(K)=W_i\otimes_{\cO} V(K)$ 
defined in \ref{subsec:The scheme H1xH2}. 
In particular we have 
\begin{align*}
\rho_i(g)\circ H^0(F^*_i)&
=H^0(F^*_i)\circ\rho'_i(g)\\
\rho_i(g)=\rho'_i(g)&=(\id_{W_i}\otimes U(K))(g).
\end{align*}

Note that Schur's lemma for $U(K)$ is true over any $\cO$-algebra. 
See \cite[Remark~7.15]{Nakamura99}.
Hence it follows from irreducibility of $U(K)$ 
that $H^0(F^*_i)=h^*_i\otimes
\op{id}_{V(K)}$ for some $h^*_i\in \GL(W_i)$. Let $\sigma(h^*_i)$ be the
transformation of $\bP(W_i(K))$ induced from $h^*_i\otimes
\op{id}_{V(K)}$. Then 
$\iota'_i\circ f_i=\sigma(h^*_i)\circ\iota_i$. 
This proves (2). (1) is clear from the proof of (2).\qed\end{pf}
%%%%%%%%%%%%%%%%%
%%%%%%%%%%%%%%%%%
\subsection{The geometric and categorical quotient}
\label{subsec:geometric and categorical quotient} 
Let $G$ be a group scheme, $X$  a scheme and $\sigma:G\times X\to X$ the
action. 
We say that 
the action $\sigma$ on $X$ is {\it proper} 
if the morphism $\Psi:=(\sigma,p_2) : G\times X\to X\times X$
is proper. We consider the following conditions:
\begin{enumerate}
\item $\phi\circ\sigma=\phi\circ p_2$,
\item $\phi$ is surjective and the image of $\Psi$ is $X\times_Y X$,
\item $\phi$ is submersive, that is, $U$ is open in $Y$ if and only if 
$\phi^{-1}(U)$ is open in $X$,
\item a given space $Z$ and a morphism $\psi:X\to Z$ such that $\psi\circ
\sigma=\psi\circ p_2$, 
then there is a unique morphism $\chi : Y\to Z$ such that
$\psi=\chi\circ\phi$. 
\end{enumerate}

A pair $(Y,\phi)$ consisting 
of an algebraic space $Y$ and a morphism $\phi : X\to Y$ 
is called a geometric quotient (resp. a categorical quotient) of $X$ 
if the conditions 1, 2 and 3 (resp. 1 and 4) are satisfied. 

The pair $(Y,\phi)$ 
is called {\it a uniform geometric quotient} (resp. {\it a uniform
categorical quotient}) if 
for any $Y$-flat $Y'$ $(Y',\phi')$ is a geometric quotient 
 (resp. a categorical quotient) of 
$X\times_S Y'$ by $G$ where $\phi':=\phi\times_YY'$.

\begin{thm}\label{thm:uniform geometric quotient}
Let $G=\GL(W_1)\times\GL(W_2)$. Then 
\begin{enumerate}
\item The action of $G$ on $U_{g,K}$ is proper and free.
\item The action of $G$ on $U_3$ is proper and has finite stabilizer.
\item The uniform geometric and uniform categorical quotient of $U_3$ by \\
$\GL(W_1)\times\GL(W_2)$ 
exists as a separated algebraic space. 
\end{enumerate}
\end{thm}
\begin{pf}Let $k$ be a closed field.
Let $Z\in U_3(k)$ and $h\in G$. 
Suppose $h\cdot Z=Z$. 
Then $h=(h_1,h_2)$ for some $h_k\in\GL(W_k)$ and 
each $h_k$ keeps $L_k$ invariant, hence $h$  
keeps $L$ invariant. This implies that $h$ is an automorphism of 
$(Z,L)$ with $G(K)$-linearization. In particular, 
$h$ induces a linear transformation $H^0(h, L)$ of $\Gamma(Z,L)$, 
which commutes with $U(K)(g)$ for any $g\in G(K)$. 
Thus $H^0(h, L)$ on $\Gamma(Z,L)$ is a scalar matrix.  \par
We assume that $H^0(h, L)$ is the identity on $\Gamma(Z,L)$. 
We shall prove that it is an automorphism of $Z$ of finite order. 
\par
First we consider the totally degenerate case, that is, 
$Z$ is a union of normal torus embeddings. We may assume
$(Z,L)=(P_0,\cL_0)$ 
by Corollary~\ref{cor:fibres over U3 and G(Z,L)}.
Since $H^0(h, L)$ is the identity on $\Gamma(Z,L)$, $h$ 
induces the identity of $(Q_0)_{\red}$ with the notation in 
Section~\ref{sec:degenerating families}. 
By Lemma~\ref{cor:freeness of Gamma(P,L)} the linear system $\Gamma(Z,L)$ 
defines a finite morphism into $\bP(\Gamma(Z,L))$, whose image is 
$(Q_0)_{\red}$.
Since $\Gamma(Z,L)$ is 
by Theorem~\ref{thm:dimension}
the $k$-vector space consisting of finite sums of monomials $\xi_a$ 
with $a\in \Del^0(Z):=$ the set of 
Delaunay vectors of $\Del(Z)$,  
each monomial $\xi_a$ in the sums is invariant under $h$. 
\par
For a given Delaunay $g$-cell $\sigma\in\Del(Z)(0)$, 
we define $X(\sigma)$ to be the sublattice of the lattice $X\ \simeq \bZ^g$
generated by Delaunay vectors of $\sigma$ 
starting from the origin.  
Let $N(\sigma)$ be the index $[X:X(\sigma)]$ and  
$N(Z)$ the least common multiple of $N(\sigma)$
for all Delaunay $g$-cells of $\Del(Z)(0)$. 
Then $\zeta_{b,0}$ is multiplied by an $N(Z)$-th root of unity 
(depending on $b$) because 
$N(Z)b\in X(\sigma)$ if $b\in C(0,\sigma)$. 
Thus $h$ is of order at most $N(Z)$. Hence the action 
of $G$ has finite stabilizer. \par
If $(Z,L)$ is a polarized abelian variety and if $H^0(h, L)$ is the
identity, 
then $h$ is the identity because $L$ is very ample.
This proves that the $G$-action on $U_{g, K}$ is free. 
The general case where 
$A_0$ in \ref{subsec:Grothendieck stable red} is nontrivial 
follows from the totally degenerate case and the fact that 
the automorphism group 
of any polarized abelian variety $(A,L)$ is finite. 
This proves that the $G$-action on $U_3$ has finite stabilizer. 
\par
It remains to prove that the action of $G$ is proper. 
This is reduced to proving the following claim
\begin{enumerate} 
\item[] Let $R$ be a discrete valuation ring $R$ 
with fraction field $k(\eta)$, 
$S=\Spec R$. Let $\sigma : G\times U_3\to U_3$ be the action and 
$\Psi=(\sigma,p_2):G\times U_3\to U_3\times U_3$.
Then for any pair $(\phi, \psi_{\eta})$ 
consisting of a morphism $\phi : S\to U_3\times U_3$ and a morphism 
$\psi_{\eta} : \Spec k(\eta)\to G\times U_3$ 
such that $\psi_{\eta}\circ\Psi=\phi\otimes_R k(\eta)$, 
there is a morphism $\psi:S\to G\times U_3$ such that 
$\psi\circ\Psi=\phi$ and $\psi
\otimes_R k(\eta)=\psi_{\eta}$.
\end{enumerate}

This is again reduced to proving the following claim:
\begin{enumerate} 
\item[] Suppose we are given an $S$-TSQAS 
$(Z,\phi_Z,\rho_Z)_{\RIG}$ with rigid $G(K)$-structure and 
a modified Mumford family 
$(P,\phi_P,\rho_P)_{\RIG}$ 
with rigid $G(K)$-structure over $S$ such that 
\begin{equation*}
(Z,\phi_Z,\rho_Z)_{\RIG}\simeq 
(P,\phi_P,\rho_P)_{\RIG}\ \op{over} \Spec k(\eta)
\end{equation*}
Then they are isomorphic over $S$. 
\end{enumerate}

In fact, the second claim follows from the proof 
of Lemma~\ref{lemma:geom fibres over U3}. 
This proves properness of the action $\Psi$.
The third assertion follows from \cite{KM97}.
\qed\end{pf}

\begin{defn}\label{reduced-coarse-moduli} 
For a contravariant functor $F$ 
from the category of algebraic spaces over $\cO$ 
to the category of sets, 
a reduced algebraic space  $W$ over $\cO$ 
and a morphism $f$ from the functor $F$ to the functor 
$h_W$ represented by $W$ is called 
{\em a reduced-coarse-moduli algebraic space} over $\cO$ of $F$ 
or we say that $F$ is 
reductively and coarsely represented over $\cO$ by $W$ 
if the following conditions are satisfied:
\begin{enumerate}
\item[(a)]   
$f(\Spec k) : F(\Spec k)\to h_W(\Spec k)$ is a bijection 
for any geometric point $\Spec k$ over $\Spec\cO$,
\item[(b)] For any reduced algebraic space $V$ over $\cO$ 
and any morphism  $g:F\to h_V$, 
there is a unique morphism $\chi : h_W\to h_V$ such that $g=\chi\circ f$.
\end{enumerate}
\end{defn}
\begin{lemma}\label{lemma:fine moduli AgK}
Assume $e_{\min}(K)\geq 3$. 
Let $A^{toric}_{g,K}$ 
be the uniform geometric quotient of $U_{g,K}$ by 
$\GL(N_1)\times\GL(N_2)$. 
Then $A^{toric}_{g,K}$ is isomorphic to the fine moduli $\cO$-scheme 
$A_{g,K}$ of $K$-symplectic abelian varieties 
in \cite{Nakamura99}.
\end{lemma}
\begin{pf}We choose and fix an
pair $d_i$ of coprime positive integers 
such that $d_i\equiv 1\mod N$ and $d_i\geq 2g+1$. Let $Y=X\times_H U_{g,
K}$ under the notation 
of \ref{subsec:The scheme U1}. Then $Y$ is  
$U_{g, K}$-flat  with fibres abelian 
varieties with level $G(K)$-structure. Since any fibre of $Y$ in the same 
$G(W_1)\times G(W_2)$-orbit
determines a unique abelian variety with rigid $G(K)$-structure 
by Lemma~\ref{lemma:existence of rigid structure}, 
we have a $G(W_1)\times G(W_2)$-invariant morphism 
$\eta: U_3\to A_{g,K}$, which induces a morphism
$\bar\eta:A^{toric}_{g,K}\to A_{g,K}$. \par
Since $A_{g,K}$ is the fine moduli, 
there is a universal family $(Z_A,\cL_A)$ of $K$-symplectic
abelian varieties with rigid $G(K)$-structure over $A_{g,K}$. 
Let $\pi_A :Z_A\to A_{g,K}$ be the natural morphism. Then 
$(\pi_A)_*(d_i\cL_A)$ is a locally free $O_{A_{g,K}}$-module. 
It is a $G(K)$-module of weight one because $d_i\equiv 1\mod N$. 
By Remark~\ref{rem:weight one, Mumford theorem}
there is a finite locally free 
$O_{A_{g,K}}$-module $W_i$ such that 
$(\pi_A)_*(d_i\cL_A)=W_i\otimes_{\cO}V(K)$ as $G(K)$-modules. 
Choosing a local trivialization of $W_i$
we have a local morphism $\eta_i : A_{g,K}\to U_{g,K}$ so that 
the composite of $\eta_i$ and 
the natural morphism of $U_{g,K}$ to $A^{toric}_{g,K}$ defines a morphism 
from $A_{g,K}$ to $A^{toric}_{g,K}$, which is the inverse of $\bar\eta$. 
This proves that $\bar\eta$ is an isomorphism. 
\qed\end{pf}

\begin{thm} Let\/  $K$\/  be 
a finite symplectic abelian group with $e_{\min}(K)\geq 3$ and 
$N=e_{\max}(K)$. The functor\/  
$\cal{SQ}^{toric}_{g,K}$\/  
is reductively and coarsely represented over $\cO_N$ by 
a complete  reduced separated algebraic space $SQ^{toric}_{g,K}$.
\end{thm}
\begin{pf} We choose and fix any pair $d_1$ and $d_2$ as before. 
Let $SQ^{toric}_{g,K}$ be the uniform geometric and uniform 
categorical quotient of $U_3$
by $\GL(W_1)\times\GL(W_2)$.  It is a complete separated 
algebraic space by Theorems~\ref{thm:summary} and \ref{thm:uniform
geometric
quotient}. Since $U_3$ is reduced, 
so is $SQ^{toric}_{g,K}$. The rest is immediate.
\qed
\end{pf}

%%%%%%%%%%%%%%%%%%%%%%%%%%%%%%%%%%%%%%%%%%%%%%%%%%%%% 
%%%%%%%%%
%%%%%%%%%
%%%%%%%%%%%%%%%%%%%%%%%%%%%%%%%%%%%%%%%%%%%%%%%%%%%%%%%%%%%%%%%%%%%%%% 

%\providecommand{\bysame}{\leavevmode\hbox to3em{\hrulefill}\thinspace} 

\end{document}